\documentclass{amsart}
\usepackage{amsfonts}
\usepackage[alphabetic]{amsrefs} 
\usepackage{color, hyperref}

\newtheorem{theorem}{Theorem}[section]
\newtheorem{lemma}[theorem]{Lemma}
\theoremstyle{definition}

\theoremstyle{remark}
\newtheorem{remark}[theorem]{Remark}
\numberwithin{equation}{section}
\theoremstyle{plain}

\newtheorem{corollary}[theorem]{Corollary}

\newtheorem{proposition}[theorem]{Proposition}

\newcommand{\R}{\mathbb{R}}
\newcommand{\calO}{\mathcal{O}}

\DeclareMathOperator{\Ric}{Ric}

\DeclareMathOperator{\Hess}{Hess}

\DeclareMathOperator{\Lie}{\mathcal{L}}
\DeclareMathOperator{\ric}{Ric}

\begin{document}
\title {Extended Solitons of the Ambient Obstruction Flow}

\author{Erin Griffin}
\address{201 Lunt Hall\\
Dept. of Math, Northwestern University\\
Evanston, IL 60208.}
\email{griffine@northwestern.edu}
\urladdr{https://sites.google.com/view/erin-griffin-math}

\author{Rahul Poddar}
\address{Harish-Chandra Research Institute\\
Prayagraj, 211019, UP, India.}
\email{rahul27poddar@gmail.com}
\urladdr{https://www.hri.res.in/people/Mathematics/rahulpoddar}

\author{Ramesh Sharma}
\address{321 Maxcy Hall\\
Dept. of Math, University of New Haven\\
West Haven, CT 06516.}
\email{rsharma@newhaven.edu}
\urladdr{https://www.newhaven.edu/arts-sciences/undergraduate-programs/mathematics/faculty.php}

\author{William Wylie}
\address{215 Carnegie Building\\
Dept. of Math, Syracuse University\\
Syracuse, NY, 13244.}
\email{wwylie@syr.edu}
\urladdr{https://wwylie.expressions.syr.edu}


\begin{abstract}
In this paper we expand on the work of the first author on ambient obstruction solitons, which are self-similar solutions to the ambient obstruction flow. Our main result is to show that any closed ambient obstruction soliton is ambient obstruction flat and has constant scalar curvature.  We show, in fact, that the first part of this result is true for a more general extended soliton equation where we allow an arbitrary conformal factor to be added to the equation.  We discuss how this implies that, on a compact manifold, the ambient obstruction flow has no fixed points up to conformal diffeomorphism other than ambient obstruction flat metrics. These results are the consequence of a general integral inequality that can be applied to the solitons to any geometric flow. Additionally, we use these results to obtain a generalization of the Bourguignon-Ezin identity on a closed Riemannian manifold and study the converse problem on a closed extended $q$-soliton. We also study this extended equation further in the case of homogeneous and product metrics. 
\end{abstract}
\maketitle

\section{Introduction}

The ambient obstruction flow was introduced by Bahuaud and Helliwell in  \cites{BHexistence, BHuniqueness}.  This geometric flow is affiliated with Fefferman and Graham's ambient obstruction tensor, $\calO$, and stationary compact metrics for the flow are ambient obstruction flat \cite{FG}. The first author studied the solitons of this geometric flow on compact and homogeneous manifolds in \cite{Griffin} and extended these results with Cunha to non-compact manifolds in \cite{CunhaGriffin}. In this paper we consider a natural extension of ambient obstruction solitons we call extended ambient obstruction solitons.

Ambient obstruction solitons are the self-similar solutions to the ambient obstruction flow. These are triples, $(M, g, X)$, of a Riemannian manifold $(M,g)$ of even dimension $n\geq 4$, and a vector field $X$ that satisfies the equation
\[ \frac{1}{2}\Lie_X g = \frac{1}{2} \Big[\calO_n +c_n (-1)^{\frac{n}{2}}\left(\Delta^{\frac{n}{2}-1} S\right)g \Big] + \lambda g \]
where $c_n$ is a constant depending on the dimension of $M$, $S$ is the scalar curvature, $\Lie$ denotes the Lie Derivative, and $\lambda$ is an arbitrary constant. A soliton is called gradient  when the vector field is the gradient of some function, $f$, on $M$. In this case $\frac{1}{2} \Lie_X g= \Hess f$.

The explicit formula for the ambient obstruction tensor is quite complicated.  (We will state it in the next section, \eqref{AOT Eqn}, where we discuss background on the ambient obstruction flow.)  On the other hand, the ambient obstruction tensor has very nice properties: it is trace free, divergence free, conformally invariant of weight $2-n$, and invariant under isometries \cites{FG, FG2, GH}. Moreover, any even dimensional Einstein metric is ambient obstruction flat and  $\calO$ is the gradient of the Q-curvature functional.  In dimension $4$, the ambient obstruction tensor is the Bach tensor, $B$, which is also the gradient of the $L^2$-norm of the Weyl curvature. This indicates that ambient obstruction flat metrics are interesting candidates for canonical metrics and the ambient obstruction flow is an interesting higher-order curvature flow whose study can shed light on the study of the ambient obstruction tensor.

 We draw the reader's attention to the scalar curvature term in the soliton equation, $c_n (-1)^{\frac{n}{2}}\left(\Delta^{\frac{n}{2}-1} S\right)$, which is conformal to the metric with conformal factor a power of the Laplacian of the scalar curvature.  Without this term the ambient obstruction flow is degenerate, \cite{BGIM}. So this term is added to counteract the conformal invariance of the tensor. The addition of this term, however, adds difficulty to the analysis of solitons.
 
 In previous work, the first author studied the case where the soliton has constant scalar curvature by focusing on homogeneous manifolds, wherein the conformal term vanishes.  She proves that compact ambient obstruction solitons with constant scalar curvature are $\calO$-flat \cite{Griffin}.  In \cite{CunhaGriffin}, Cunha and the first author extend this notion to non-compact manifolds, finding similar results when additional conditions at infinity are imposed on $X$. In this paper we extend some of this analysis to the general case, in fact we can handle any conformal factor added to the equation.

 We define an {\bf extended ambient obstruction soliton}, to be  $(M, g, X)$ satisfying 
\begin{equation}
\frac{1}{2}\Lie_X g = \frac{1}{2}\calO + \varphi g 
\end{equation}
for some  $\varphi: M \to \R$. Note that an ambient obstruction soliton is an extended ambient obstruction soliton with $\varphi = \lambda +  \frac{1}{2} c_n (-1)^{\frac{n}{2}}\left(\Delta^{\frac{n}{2}-1} S\right)$.  We have used the nomenclature ``extended" soliton following \cite{KDT} where the concept is introduced for the Cotton flow. For the Ricci flow, a similar notion is also called a Ricci almost soliton, \cite{PRRS}. 

Our main result is that there are, in fact, no extended ambient obstruction solitons on compact manifolds aside from the ambient obstruction flat ones.

\begin{theorem}\label{Thm: GAOS Compact}
Any extended ambient obstruction soliton on a closed manifold is ambient obstruction flat. 
\end{theorem}

\begin{remark} In fact, our argument works for any flow where the trace-free part of the tensor defining the flow is divergence free and follows from an integral identity that is true for extended solitons to any geometric flow, see Theorem \ref{Thm:compact} in Section 3.
\end{remark}

We also note that when we consider a soliton which is ambient obstruction flat we get the following rigidity.

\begin{theorem} \label{Thm:AOSCompact}
Any ambient obstruction soliton on a closed manifold is ambient obstruction flat, has constant scalar curvature, and its soliton vector field $X$ is Killing.
\end{theorem}

In the non-compact case there are constant scalar curvature, in fact homogeneous, ambient obstruction solitons that are not ambient obstruction flat.  For example, Ho showed there are gradient product metrics examples on $\mathbb{R}^2 \times S^2$ and  $\mathbb{R}^2 \times H^2$ \cite{Ho}, the first author found a gradient product soliton on $\mathbb{R} \times SU(2)$ \cite{Griffin},  and Thompson finds a non-gradient example on $N^4$,  the simply connected $4$-dimensional nilpotent Lie group \cite{Thompson}. 

Our next result shows that, in the homogeneous case, extended ambient obstruction solitons are solitons. 

\begin{theorem}\label{Thm: GAOS Hom}
Any extended ambient obstruction soliton on a homogeneous space is an ambient obstruction soliton. 
\end{theorem}

We note that this ``extended" notion of soliton also has a natural geometric interpretation in terms of the original ambient obstruction flow.  Namely, just as solitons describe the solutions to the flow that are fixed points modulo diffeomorphism and re-scaling, any flow that is a fixed point modulo diffeomorphism and conformal change must be a extended ambient obstruction soliton.   We discuss this further in Section \ref{Conf Diff}.

 When $n=4$ the ambient obstruction tensor is the Bach tensor, $B$.  Consequently, in dimension $4$ we call ambient obstruction solitons Bach solitons.  The soliton equation in this case is 
\[ \frac{1}{2} \Lie_X g = \frac{1}{2} \left( B + \frac{1}{12} \Delta S g\right) + \lambda g.  \]
The first author also studied homogeneous gradient Bach solitons in dimension $4$ in \cite{Griffin}. By a result of the fourth author and Petersen, any such metric that is not Bach flat must be a product of a Euclidean space and another factor, $N$ \cite{PW}.   It follows from the results in \cite{Griffin} that the only locally homogeneous gradient Bach solitons where the factor $N$ is compact are quotients of a unique (up to scaling) metric on $\mathbb{R} \times SU(2)$  and metrics  found by Ho of the form $\mathbb{R}^2 \times N^2$ where $N^2$ is a space of constant curvature.  Motivated by this, we consider the case where we have a product gradient Bach soliton with a compact factor that is not necessarily homogeneous. We do not obtain a classification but do obtain the following  structure result. 

\begin{theorem}
Let $(M,g,f)$ be a non-Bach flat gradient Bach soliton on a product metric 
\begin{enumerate}
    \item If $M = S^1 \times N^3$ then $\lambda>0$ and $f$ is a function on $N$ only. 
    \item If $M = \mathbb{R} \times N^3$ where $N^3$ is compact, then $\lambda <0$, $f= 2 \lambda t^2 +  a t + b$, and  $N$ has constant scalar curvature and constant norm of the Ricci tensor. 
    \item If $M = K^2 \times L^2$ where $K^2$ is compact then $\lambda \geq 0$,  $g_K$ is constant curvature, and $f$ is a function on $L$ only. 
\end{enumerate}
\end{theorem}

\begin{remark} \label{HoRemark} The examples of Ho show that case (3) can appear and the first author's example on $\mathbb{R} \times SU(2)$ shows that case (2) can appear \cite{Ho, Griffin}. We do not know if there are examples in case (1).  In cases (2) and (3) there are local solutions that are not isometric to Ho or the first author's examples, but we do not know whether Ho and the first author's are the only complete solutions.
\end{remark}

The paper is organized as follows.  In the next section we give more background detail on the ambient obstruction tensor and extended solitons.  In section 3 we prove the general theorems that imply Theorems \ref{Thm: GAOS Compact} and \ref{Thm: GAOS Hom}. These are also applied to obtain a generalization of the Bourguignon-Ezin identity on a closed Riemannian manifold and study the converse problem on an extended $q$-soliton. In section 4 we discuss product metrics, where, in addition to the result stated above for gradient solitons, we also obtain some partial results for non-gradient solitons and extended solitons.

\section{Background}

The ambient obstruction tensor, $\calO$, defined for even dimension $n\geq 4$, is an $n$th order symmetric $2$-tensor. It is the obstruction to a manifold having a Fefferman-Graham ambient metric associated with the given conformal structure \cite{GH}. Thus, an $\calO$-flat metric allows the manifold to have a smooth formal power series solution for the ambient metric associated to the given conformal structure.

The ambient obstruction tensor is given explicitly by the equation
\begin{equation}\label{AOT Eqn}
\begin{gathered}  
\calO_n = \frac{1}{(-2)^{\frac{n}{2}-2} \left(\frac{n}{2}-2 \right)! } \left( \Delta^{\frac{n}{2}-1} P - \frac{1}{2(n-1)} \Delta^{\frac{n}{2}-2} \nabla^2 S \right) + T_{n-1} \\
 P = \frac{1}{n-2} \left( \ric - \frac{1}{2(n-1)} S_g \right) \end{gathered}\end{equation}
where $P$ is the Schouten tensor and $T_{n-1}$ a polynomial natural tensor of order $n-1$ \cite{FG, GH, Lopez}. Importantly, it is well established that $\calO$ is trace-free, divergence-free, and conformally invariant of weight $2-n$. 

The ambient obstruction flow is given by the equation, 
\begin{equation} \label{AOF}
\begin{cases} 
\partial_t g = \calO_n +c_n (-1)^{\frac{n}{2}}\left(\Delta^{\frac{n}{2}-1} S\right)g \\
g(0)=h,
\end{cases}\end{equation}
where $h$ is a smooth metric of even dimension $n\geq 4$, $S$ denotes the scalar curvature,  and $c_n$ is a dimension-dependent constant. In dimension $4$, we will call flow the Bach flow and equation becomes 
\begin{equation*}
\partial_t g = B + \frac{1}{12} \left(\Delta S\right)g.
\end{equation*}

Bahuaud and Helliwell define the ambient obstruction flow, and show its short time existence and uniqueness on compact manifolds in \cite{BHexistence, BHuniqueness}, respectively.  An alternate proof to the uniqueness of solutions to (\ref{AOF}) was also given by Kotschwar \cite{Kotschwar}. Recently Bahuaud, Guenther, Isenberg, and Mazzeo showed existence and uniqueness on complete manifolds with bounded geometry \cite{BGIM}. Lopez also shows that, under the ambient obstruction flow, the curvature blows up at a finite singular time and studies singularity models and non-singular solutions \cite{Lopez}.

Our proofs do not rely on the explicit formula for $\mathcal{O}$, but on its property of being both divergence and trace free. Consequently, they apply more broadly to other flows with this property. Explicitly,  if $q$ is a symmetric two-tensor, then a $q$-soliton is some $(M, g, X)$ satisfying an equation of the form $\frac{1}{2}\Lie_X g = \frac{1}{2} q + \lambda g$ for some constant $\lambda$. If we replace the constant $\lambda$ with a function we will also call the corresponding a equation the {\bf extended $q$-soliton} equation:
\begin{equation}\label{GQS} \frac{1}{2}\Lie_X g = \frac{1}{2}q + \varphi g. \end{equation}

\subsection{Fixed points up to conformal diffeomorphisms and extended solitons}\label{Conf Diff}

In this section we give an interpretation of what the extended soliton equation in terms of the corresponding geometric flow.  The basic idea is that while the usual soliton equation is the infinitesimal version of a flow that changes only by homotheties, the extended equation is the infinitesimal version of a flow that changes only by conformal diffeomorphisms.

The idea of modifying the constant in the soliton equation to be a function is well known.  In the case of the Ricci flow, what we have called an extended soliton is called an ``almost Ricci soliton''.  In that case there are a number of rigidity results known, however, we are not aware of a place in the literature that discusses the interpretation in terms of conformal diffeomorphisms.

A diffeomorphism between Riemannian manifolds $F: (M_1, g_1) \to (M_2, g_2)$ is called \emph{conformal} if it preserves angles, equivalently, $F^*(g_2) = v g_1$ for some positive function $v$.   Compositions and inverses of conformal diffeomorphisms are conformal diffeomorphisms. Therefore, conformal diffeomorphism is an equivalence relation on the space of metrics. The conformal diffeomorphism class of a metric is the set of metrics that are related to it by a conformal diffeomorphism. 

We say a flow (that is, a smooth one-parameter family of metrics, $g(t)$) is a fixed point up to conformal diffeomorphism if $g(t)$ is in the same conformal diffeomorphism class for each $t$. Then we can write  
\begin{eqnarray}
g(t) =v_tF_t^*(g(0)). \label{CDFeqn}
\end{eqnarray} 
for some functions $v_t: M \to R$ and one-parameter family of diffeomorphisms $F_t: M \to M$ such that $v_0 = 1$ and $F_0 = id$. We will assume the functions $v_t$ and $F_t$ are differentiable. 

 Note that a fixed point up to conformal diffeomorphism is a generalization of a soliton metric, which is a flow by homotheties, where $v_t$ depends only on $t$. Understanding solitons for geometric flows is important because they act as the geometric fixed points, as the underlying geometry of a soliton is not changed by the flow, and thus can not be ``improved" by the flow.   Similarly, for a fixed point up to conformal diffeomorphism,  the conformal properties of the metric will not improve under the flow. Thus fixed points up to conformal diffeomorphism seem to be of natural interest to the application of geometric flows to conformal geometry.

We now prove our main observation about conformal diffeomorphisms  and extended solitons.  We say that $g(t)$ is a solution to the $q$-flow if $\frac{\partial g}{\partial t} = q(g(t))$, where $q$ is a map from metrics to symmetric $2$-tensors. Differentiating both sides of equation (\ref{CDFeqn}) we obtain the following proposition. 

\begin{proposition} \label{prop:CDGS}
    If $g(t)$ is a solution to the $q$-flow that is fixed up to conformal diffeomorphism then it is an extended $q(g(t))$-soliton. 
\end{proposition}

\begin{proof}
Differentiating both sides of (\ref{CDFeqn})
 with respect to $t$ at $t=0$ gives
\begin{eqnarray*} 
\frac{\partial g}{\partial t}\Big\vert_{t=0} = \frac{\partial v}{\partial t}\Big\vert_{t=0} g(0) + \Lie_X g(0)
\end{eqnarray*}
So if $g$ is a solution to $q$ flow then we have
\begin{eqnarray*} 
q(g(0)) = \frac{\partial v}{\partial t}\Big\vert_{t=0} g(0) + \Lie_X g(0)
\end{eqnarray*}
which is equivalent to the extended $q$-soliton equation.
\end{proof}

We can also prove a converse statement if $q$ is conformally invariant. 

\begin{proposition} \label{Prop:CDGS converse} Suppose that $g_0$ is an extended $q(g_0)$-soliton for some diffeomorphism invariant $q$ with the property that $q(vg) = vq(g)$ for all $v$ and $g$.  Then there is a solution to the $q$-flow, $g(t)$ with $g(0) = g_0$ that is fixed up to conformal diffeomorphism. 
\end{proposition}

\begin{proof}
If $g(t) = v_t F^*_t(g_0)$ then by the product rule
\[ \frac{dg}{dt}|_{t=t_0} = \left(\frac{\partial v}{\partial t}|_{t=t_0} \right) F_{t_0}^*(g_0) + v_{t_0} \Lie_{Y(t_0)}( F_{t_0}^* g_0)\]
where $Y(t_0) = F_{t_0}^* ( \frac{\partial}{\partial t} |_{t=t_0} F_t)$. On the other hand, since $q(vg)= vq(g)$ and $q$ is diffeomorphism invariant 
\begin{align*}
    q(g(t)) &= q(v_t F_t^*(g_0)) \\
    &= v_t F_t^*(q(g_0)) \\
    &= v_t F_t^*( \Lie_X g_0 + \varphi g_0) \\
    &= v_t \Lie_{F_t^*(X)} (F_t^*(g_0)) + v_t F_t^*(\varphi) F_t^*(g_0)
\end{align*}
We thus see that we have a flow that is fixed by conformal diffeomorphism if we set $\frac{\partial}{\partial t} |_{t=t_0} F_t = X$ and then have $v$ solve $\frac{\partial v}{\partial t} = v F_t^*(\varphi)$.

\end{proof}

We note that Proposition \ref{Prop:CDGS converse} does not seem to be true unless $q$ is exactly conformally invariant. In particular,  it is not clear if every extended ambient obstruction soliton gives rise to a fixed point up to conformal diffeomorphism. The extended soliton equation is the infinitesimal obstruction to the flow being fixed point up to conformal diffeomorphism that comes from differentiation, but it is not clear whether the extended soliton equation can always be ``integrated" to give a corresponding flow.   This is in contrast to the regular soliton equation, where the infinitesimal soliton equation will give rise to a  flow that is a fixed point under homotheties as long as $q$ satisfies a re-scaling property.  But this is exactly because in the case of homotheties, the re-scaling property will give an ODE for the re-scaling function that can always be solved, but in the  conformal case, one gets a PDE that $v_t$ must solve which does not seem always possible. 

On the other hand,  non-existence of extended solitons on a given manifold obstructs the existence of a flow that is fixed under conformal diffeomorphism. Therefore, Theorem \ref{Thm: GAOS Compact} shows that the ambient obstruction flow on a compact manifold has no fixed point under conformal diffeomorphism other than the ambient obstruction flat metrics.

\section{Proof of main theorems}

In this section we prove the results that imply our main theorems for extended ambient obstruction solitons, Theorems \ref{Thm: GAOS Compact}, \ref{Thm:AOSCompact} and \ref{Thm: GAOS Hom}.  In fact we prove results in the general setting of extended solitons to any diffeomorphism invariant flow, implying the results in the specific case of the ambient obstruction flow. 

First we give two integral formulas that imply Theorem \ref{Thm: GAOS Compact}. 
\begin{theorem} \label{Thm:compact}
Let $(M,g,X)$ be an extended $q$-soliton on a closed Riemannian manifold, 
\[ \frac{1}{2} \Lie_X g = \frac{1}{2}q + \varphi g\]
where $\varphi$ is some function on $M$.  Let $\mathring{q} = q - \frac{ \mathrm{tr}(q)}{n}g$, the trace free part of $q$. Then 
\begin{align} \label{eqn:integral1} \int_M (\varphi \mathrm{tr}(q) + \frac{(\mathrm{tr}(q))^2}{2n} + (\mathrm{div} q)(X)) dvol_g &= -  \frac{1}{2}\int_M ||\mathring{q}||^{2} dvol_g \leq 0 \\ \label{eqn:integral2} \int_M (\mathrm{div} \mathring{q})(X) dvol_g &= -\frac{1}{2} \int_M || \mathring{\Lie_X g}||^2dvol_g \leq 0. \end{align}
In particular, the integral on the left hand side of either equation is  equal to $0$ if and only if $X$ is a conformal field and $q= ( \frac{\mathrm{div}X}{n} - \varphi)g$.  
\end{theorem} 

Before proving the theorem, we point out the main corollary we are after. 

\begin{corollary} \label{Cor:compact} Let $(M,g,X)$ be an extended $q$-soliton on a closed Riemannian manifold such that $q$ is divergence free. 
\begin{enumerate}
\item If  $q$ has constant trace, then $X$ is a conformal field and $q= ( \frac{\mathrm{div}X}{n} - \varphi)g$.
\item If  $\mathrm{div}(X) = 0$, then $X$ is Killing and $q$ is a multiple of the metric.
\item  If $q$ is trace free or $\varphi=0$, then $X$ is Killing and $q=0$. 
\end{enumerate}
\end{corollary} 

\begin{proof}[Proof of Corollary \ref{Cor:compact}]
If $q$ is divergence free and constant trace, then $\mathrm{div}(\mathring{q})=0$ so the conclusion follows from  (\ref{eqn:integral2}).  If $q$ is divergence free and trace free then (\ref{eqn:integral1}) implies that $||\mathring{q}||=0$ which implies $q=0$ since $\mathrm{tr}(q)=0$.  If $q$ is divergence free and $\varphi=0$ then  (\ref{eqn:integral1}) implies that $\mathrm{tr}(q) = 0$, so we also have $q=0$ in this case.  Finally, if $\mathrm{div}(X) =0$ then from the trace of the extended $q$-soliton equation we have that 
\[ \mathrm{div}(X) = \frac{1}{2} \mathrm{tr}(q) + n \varphi\]
So $\mathrm{div}(X) = 0$ is equivalent to $\mathrm{tr}(q) = -2n \varphi$.  But this implies that $\varphi \mathrm{tr}(q) + \frac{(\mathrm{tr}(q))^2}{2n} = 0$.  So  (\ref{eqn:integral1}) implies that $\mathring q = 0$.  Since $\mathrm{div}(X) = 0$ this implies that $X$ is Killing. 

\end{proof}

We note that since the ambient obstruction tensor is divergence and trace free, Corollary \ref{Cor:compact} implies Theorem \ref{Thm: GAOS Compact}. Now we prove Theorem \ref{Thm:compact}.

\begin{proof}[Proof of Theorem \ref{Thm:compact}]
    
    Both equations (\ref{eqn:integral1}) and (\ref{eqn:integral2}) follow from the following identity: For a symmetric $(0,2)$-tensor field $T$ and an arbitrary smooth vector field $X$ on $M$, we have the identity \cite[Lemma 3.5]{Griffin}:
\begin{equation} \label{3.5}
        \langle \Lie_X g, T \rangle = 2\mathrm{div}(i_X T) - 2(\mathrm{div} T)(X).
\end{equation}

To establish (\ref{eqn:integral1}), we let $T=q$ and substitute $\Lie_X g = \frac{1}{2} q + \phi g$ into the left hand side to obtain
\begin{equation*}
    ||q||^{2} + 2\varphi \mathrm{tr}(q) = 2\mathrm{div}(i_X q) - 2(\mathrm{div} q)(X).
\end{equation*}
Integrating it over compact $M$ and using $||\mathring{q}||^{2} = ||q||^{2} - \frac{(\mathrm{tr}(q))^2}{n}$ we have that
\begin{equation*}
     \int_M (\frac{1}{2} ||\mathring{q}||^{2} + \varphi \mathrm{tr}(q) + \frac{(\mathrm{tr}(q))^2}{2n} + (\mathrm{div} q)(X)) dvol_g = 0.
\end{equation*}

  On the other hand, if we let $T= \Lie_X g$ in (\ref{3.5}) and integrate we obtain 
    \begin{align} \label{eqn:LieIntegral} \int_M ||\Lie_X g||^2 dvol_g = -2 \int_M \mathrm{div}(\Lie_X g)(X) dvol_g 
    \end{align}
Taking the divergence of the generalized soliton equation, and using the trace   $\mathrm{div}(X) = \frac{1}{2} \mathrm{tr}(q) + n \varphi$, we also obtain 
    \begin{align}
    \mathrm{div}(\Lie_X g) &=\mathrm{div}(q) + 2 \mathrm{div}\left( \varphi g \right) \nonumber \\
    &= \mathrm{div}(q) +  2 d\varphi  \nonumber  \\
    &= \mathrm{div}(q) + \frac{2}{n} d\left(\mathrm{div} X - \frac{1}{2} \mathrm{tr}(q)  \right). \nonumber \\
    &= \mathrm{div}( \mathring{q})  + \frac{2}{n} d( \mathrm{div}(X))\label{eqn:div}
    \end{align}
    For any function $h$ we have the identity  $\mathrm{div}(hX) = dh(X) + h \mathrm{div}(X)$, applying this to $h=\mathrm{div}(X)$ we obtain
    \begin{align*}
        d(\mathrm{div}X)(X) = \mathrm{div}(\mathrm{div}(X) X) - (\mathrm{div}(X))^2.
    \end{align*}
    Integrating both sides then gives that 
    \begin{align*}
        \int_M d(\mathrm{div}X)(X) = - \int_M (\mathrm{div}X)^2 dvol_g. 
    \end{align*}
    Putting these this together with (\ref{eqn:LieIntegral}) and (\ref{eqn:div}) yields that 
    \begin{align*}
        \int_M ||\Lie_X g||^2 dvol_g =  -2 \int_M \mathrm{div}(\mathring{q})(X) dvol_g +  \frac{4}{n} \int_M (\mathrm{div}X)^2 dvol_g. 
    \end{align*}
    So that
    \begin{align*}
      2\int_M \mathrm{div}(\mathring{q})(X) dvol_g  = -  \int_M \left( ||\Lie_X g||^2- \frac{4}{n}(\mathrm{div}X)^2\right) dvol_g.
    \end{align*}
    Since the trace of $\Lie_X g$ is $2 \mathrm{div}(X)$, this gives (\ref{eqn:integral2}).  
    
    Now if the left hand side of (\ref{eqn:integral1}) is zero, then $\mathring q =0$, so $q = \frac{\mathrm{tr}{q}}{n} g = 2\left( \frac{\mathrm{div} X}{n} - \phi \right) g.$ If the left hand side of (\ref{eqn:integral2}) is zero then we have that $X$ is a conformal field.  Then we have $\Lie_Xg = \frac{2}{n} \mathrm{div}(X) g$, so plugging back into the extended soliton equation we get that $q= 2( \frac{\mathrm{div}X}{n} - \varphi)g$, in this case as well. 

\end{proof}

\begin{remark}   In the case of almost Ricci solitons, we have $q = - 2\Ric$. Recall that by the second Bianchi identity, $\mathrm{div}(\Ric) = \frac{1}{2} d S$, so in this case  
\[ \mathrm{div}\mathring{q} = \frac{2-n}{n} dS.  \] Consequently,  Theorem \ref{Thm:compact} implies that any compact almost Ricci soliton with constant scalar curvature is Einstein.  In fact, Theorem \ref{Thm:compact} was proven in \cite{BBR} for Ricci almost solitons.  

Further, we can see that
\begin{equation*}
    \int_M (\mathrm{div}(Ric))X dvol_g = \frac{1}{2}\int_M (\Lie_X S) dvol_g = 0,
\end{equation*}
in view of $\mathrm{div}(SX)=\Lie_X S+S\mathrm{div} X$, compact $M$ and $\mathrm{div} X = 0$. Thus, part (2) of Corollary \ref{Cor:compact} implies the following result: ``Any compact almost Ricci soliton with divergence free $X$ is Einstein and $X$ is Killing" proved in \cite{Sharma} for an almost Ricci soliton. 
\end{remark}

\begin{proof}[Proof of Theorem  \ref{Thm:AOSCompact}] By Theorem \ref{Thm: GAOS Compact} we have $\calO=0$. So the 
ambient obstruction soliton equation defined in the beginning of Section 1 reduces to
\begin{equation}\label{3.4}
    \frac{1}{2}\Lie_X g = \lambda g + \frac{1}{2}c_n (-1)^\frac{n}{2}(\Delta^{\frac{n}{2}-1}S)g.
\end{equation}
Taking the trace:
\begin{equation*}
    \mathrm{div} X = n\lambda + \frac{n}{2}c_n(-1)^\frac{n}{2}\Delta(\Delta^{\frac{n}{2}-2}S).
\end{equation*}
Integrating over compact $M$ immediately gives $\int_M \lambda =0$, so $\lambda=0$. Consequently, equation (\ref{3.4}) becomes the conformal equation
\begin{equation*}
    \Lie_X g = 2\sigma g,
\end{equation*}
where $\sigma = \frac{1}{n}\mathrm{div} X=\frac{c_n}{2}(-1)^\frac{n}{2}(\Delta^{\frac{n}{2}-1} S)$. Using this in the conformal integrability condition of Yano, \cite{Yano}: $\Lie_X S=-2\sigma S -2(n-1)\Delta \sigma$ ($\Delta \sigma=\mathrm{div} (\nabla \sigma)$) and also using the Bourguignon-Ezin identity \cite{Bourgignon-Ezin}: $\int_M (\Lie_X S) dvol_g=0$, provides $\int_M S \sigma dvol_g =0$, i.e.,
\begin{equation}\label{1.13}
    \int_M S( \Delta^{\frac{n}{2}-1}S)dvol_g = 0.
\end{equation}
At this point, we use the self-adjointness of $\Delta$ successively, on a compact Riemannian manifold, moving $\Delta$'s a sufficient number of times to the left factor of the above equation in order to get
\begin{equation*}
\int_M (\Delta^{p}S)^2 dvol_g = 0
\end{equation*}
or
\begin{equation*}
  \int_M (\Delta^{p}S) (\Delta^{p+1}S) dvol_g = 0  
\end{equation*}
according as $\frac{n}{2}-1 = 2p$ or $2p+1$ for a non-negative integer $p$. In the first case, $S$ is constant by Hopf's lemma. In the second case, we compute
\begin{eqnarray*}\label{1.16}
   \Delta ((\Delta^{p}S)^{2}) &=& \nabla^{i}\nabla_{i}((\Delta^{p}S)^{2})\nonumber\\ &=& \nabla^{i}[2(\Delta^{p}S)\nabla_{i}\Delta^{p}S]\nonumber\\ &=& 2[||\nabla(\Delta^{p}S)||^{2} + (\Delta^{p}S)(\Delta^{p+1}S)].
\end{eqnarray*}
Integrating it over $M$ and using the equation of the second case, shows that $\nabla(\Delta^{p}S) = 0$, i.e., $\Delta^{p}S$ is constant, which again, by Hopf's lemma, implies that $S$ is constant. This completes the proof.
\end{proof}

Now we consider our general result in the homogeneous case that implies Theorem \ref{Thm: GAOS Hom}

\begin{theorem} \label{Thm:homogeneous} Suppose that $(M,g)$ is a homogeneous extended $q$-soliton  \[ \frac{1}{2} \Lie_X g = \frac{1}{2}q + \varphi g\] for some isometry invariant $2$-tensor $q$, then either $M$ is locally conformally flat or $(M,g)$ is a $q$-soliton. \end{theorem}

We note that locally conformally flat homogeneous spaces are classified by Takagi \cite{Takagi}. To see that Theorem \ref{Thm:homogeneous} implies Theorem \ref{Thm: GAOS Hom} note that by the conformal invariance of $\mathcal{O}$, locally conformally flat metrics are always ambient obstruction flat.  Theorem \ref{Thm:homogeneous} was also observed in the Ricci almost soliton case in \cite{HomRAS}. 
\begin{proof}
 Let $F$ be an isometry of $(M,g)$.  Then
\begin{eqnarray*}
\frac{1}{2} F^*(\Lie_X g) &=& \frac{1}{2} F^*q + F^* (\varphi g)\\
\frac{1}{2}\Lie_{(F^* (X))}g &=& \frac{1}{2} q + (F^* \varphi) g.
\end{eqnarray*}
Subtracting $\frac{1}{2}\Lie_X g = \frac{1}{2}q + \varphi g$ from this equation gives us 
\[ \frac{1}{2} L_{(F^* (X)-X)}g = (F^* (\varphi) - \varphi) g \]
If $(F^* (\varphi) - \varphi) \neq 0$ then this makes $F^*(X) - X$ a conformal field. By \cite[Proposition 2.5]{PW}  if a homogeneous space admits a non-Killing conformal field, then the space is locally conformally flat.  Therefore, if the space is not locally conformally flat, then $F^* (\varphi)= \varphi$ for all isometries $F$.  Then, since the space is homogeneous $\varphi$ must be constant. 
\end{proof}

We note that there are Ricci almost solitons with $\varphi$ non-constant that are not homogeneous and not locally conformally flat, see \cite[Corollary 1]{FFGP}.

We would like to introduce an equation that considers both the Ricci flow and ambient obstruction flow simultaneously. That is,  we assume that $q$ is a second order symmetric tensor $q$ on $M$ satisfying
\begin{equation}\label{3.7}
    \mathrm{div} q = \frac{1}{2}\nabla \mathrm{tr}(q).
\end{equation}
By the twice contracted second Bianchi identity, $\mathrm{div} Ric = \frac{1}{2}\nabla S$, so $\mathrm{Ric}$ is an example of such a tensor.  Moreover, so is any divergence free tensor with constant trace, such as the ambient obstruction tensor.  The  Bourguignon and Ezin  formula \cite{Bourgignon-Ezin} for a conformal vector field $X$ on a compact Riemannian manifold $M^n$ ($n>2$) says that  $\int_M (\Lie_X S) dvol_g = 0$. We give an extension of this result to other $q$ satisfying (\ref{3.7}).

\begin{theorem}
    Let $X$ be a conformal vector field on a compact Riemannian manifold $M^n$ ($n>2$), and $q$ a symmetric $(0,2)$ tensor such that $\mathrm{div} q = \frac{1}{2}\nabla \mathrm{tr}(q)$. Then $M$ satisfies the integral formula
    \begin{equation*}
        \int_M (\Lie_X \mathrm{tr}(q)) dvol_g = 0.
    \end{equation*}
\end{theorem}

\begin{proof}
Since $X$ is conformal, we have 
\begin{equation*}
    \Lie_X g = 2\sigma g,
\end{equation*}
where $\sigma$ is the conformal scale function. Using this and the hypothesis $\mathrm{div} q = \frac{1}{2}\nabla \mathrm{tr}(q)$, in the identity used in \cite[Lemma 3.5]{Griffin}, we observe that: 
\begin{equation*}
    \langle \Lie_X g, q \rangle = 2\mathrm{div}(i_X q) - 2(\mathrm{div} q)(X)
\end{equation*}
and then integrating over the compact $M$ provides
\begin{equation*}
    \int_M (2\sigma \mathrm{tr}(q) + \Lie_X \mathrm{tr}(q))dvol_g = 0.
\end{equation*}
Now, integrating the identity: $ \mathrm{div}(\mathrm{tr}(q)X) = \Lie_X \mathrm{tr}(q) + \mathrm{tr}(q)\mathrm{div} X$ over $M$ and noting that $\mathrm{div} X = n\sigma$ (trace of the conformal equation), we have
\begin{equation*}
    \int_M (n\sigma \mathrm{tr}(q) +\Lie_X \mathrm{tr}(q))dvol_g = 0.
\end{equation*}
Comparing the last two integral equations, and noting $n>2$, we get the integral formula of the theorem, completing the proof.

\end{proof}

\noindent
Finally, we show that the converse of the previous result holds if the compact Riemannian manifold is an extended $q$-soliton. More precisely, we prove the following result.
\begin{theorem}
   If a closed extended $q$-soliton $(M^n,g,X,q)$ satisfies the conditions:
    \begin{equation}\label{5.3}
         \mathrm{div} q = \frac{1}{2}\nabla \mathrm{tr}(q) \hspace{0.2cm}, \int_M (\Lie_X \mathrm{tr}(q)) dvol_g = 0
    \end{equation}
then $X$ is a conformal vector field. For dimension $n=2$, the second condition is redundant.
\end{theorem}

\begin{proof} Using the value of $\varphi$ from the trace of \eqref{GQS} in equation \eqref{eqn:integral1} (which is applicable here) gives
\begin{equation*}
     \int_M \left(\frac{1}{2} ||\mathring{q}||^{2} + \frac{1}{n}\mathrm{tr}(q)\mathrm{div} X + (\mathrm{div} q)(X)\right) dvol_g = 0.
\end{equation*}
Next, the integral of the identity: $\mathrm{div}(\mathrm{tr}(q)X) = \Lie_X \mathrm{tr}(q) + \mathrm{tr}(q)\mathrm{div} X$ over $M$ provides
\begin{equation*}
    \int_M \left(\Lie_X \mathrm{tr}(q) + \mathrm{tr}(q)\mathrm{div} X\right)dvol_g = 0.
\end{equation*}
The preceding two equations show that
\begin{equation*}
     \int_M \left(\frac{1}{2} ||\mathring{q}||^{2} - \frac{1}{n}\Lie_X \mathrm{tr}(q) + (\mathrm{div} q)(X)\right) dvol_g = 0.
\end{equation*}
The use of the hypothesis: $\mathrm{div} q = \frac{1}{2}\nabla \mathrm{tr}(q)$ in the above equation provides
\begin{equation}\label{5.4}
    \int_M \left(||\mathring{q}||^{2} + \frac{n-2}{n}\Lie_X \mathrm{tr}(q)\right)dvol_g = 0.
\end{equation}

Finally, using the second integral equality (\ref{5.3}) of our hypothesis in the above equation, we obtain $\mathring{q}=0$. This consequence, in conjunction with (\ref{GQS}) shows that $X$ is conformal. For $n=2$, equation (\ref{5.4}) gives the same conclusion without assuming the second integral equality (\ref{5.3}) of our hypothesis. This completes the proof.
\end{proof}

\section{Bach solitons on products} \label{bach section}

In this section we consider the Bach soliton equation on products. We begin with  a general lemma about when vector fields and functions split on products.  Let $M = N_1 \times N_2$ be a product of two Riemannian manifolds with a product metric $g_M = g_{N_1} + g_{N_2}$.  Let $X$ be a vector field on $M$.  We will say that $X$ splits if $X = X_1 + X_2$ where $X_1$ is a vector field on $N_1$ and $X_2$ is a vector field on $N_2$.  We say a function $f:M \to \mathbb{R}$ splits if $f= f_1 + f_2$ where $f_i:N_i \to \mathbb{R}$.  

We have the following Lemma about splitting of functions and vector fields in terms of the Lie derivative of the metric $g$. 

\begin{lemma} \label{Lem:Splitting} If $\mathrm{Hess} f(Y,U) = 0$ for all vectors $Y\in N_1$ and $U \in N_2$ then $f$ splits.  If $X$ is a vector field such that $(\Lie_X g)(Y,U) = 0$ then $X$ splits if and only if $\nabla_Z X \in T_p N_1$ for all $Z \in N_1$ (equivalently $\nabla_U X \in T_p N_2$ for all $U \in N_2$). 
\end{lemma}

\begin{proof}
    We prove the statement for non-gradient fields first. 
    Let $\{\frac{\partial}{\partial x_i}\}$, $\{ \frac{\partial}{\partial y_{a}}\}$ be coordinate vector fields on $N_1$ and $N_2$ respectively that are orthonormal and have zero covariant derivative  at a point $(p,q)$.  Write
    \[ X = X^i(x,y)\frac{\partial}{\partial x_i} + X^{a}(x,y) \frac{\partial}{\partial y_{a}} \]
    in a neighborhood of $(p,q)$. $X$ is splitting in the neighborhood is equivalent to  $X^i$ being a  function of $x$ only and $X^a$ being  function of $y$ only. The assumption $\nabla_Z X \in T_p N_1$ implies that $X^a$ is a function of $y$ for each $a$ since we then have  
    \begin{align*}
         0 = g\left(\nabla_{\frac{\partial}{\partial x_j}} X, \frac{\partial}{\partial y_{a}}\right) = \frac{\partial X^{a}}{\partial x_j}.
    \end{align*}
   The mixed term in the Lie Derivative is 
    \begin{align*}
        (\Lie_X g)\left(\frac{\partial}{\partial x_j},\frac{\partial}{\partial y_{a}}\right)&= g\left(\nabla_{\frac{\partial}{\partial x_j}} X, \frac{\partial}{\partial y_{a}}\right) + g\left(\frac{\partial}{\partial x_j}, \nabla_{\frac{\partial}{\partial y_{a}}} X\right) \\
        &= \frac{\partial X^{a}}{\partial x_j} +  \frac{\partial X^{j}}{\partial y_a}. \end{align*} 
    So we see that under the assumption that $0 =(\Lie_X g)\left(\frac{\partial}{\partial x_j},\frac{\partial}{\partial y_{a}}\right)$,  $X^a$ is a function of $y$ only for all $a$  is equivalent to $X^i$  being a function of $x$ only for all $i$.  So we obtain that $X$ is splits.

    In the gradient case, this fact has also been observed in \cite{PWSymmetry}, but it is also easy to see from the above, so we give the proof for completeness.  Now we have $X^i = \frac{\partial f}{\partial x_i}$, $X^a = \frac{\partial f}{\partial y_a}$, so 
    \begin{align*}
        0 = \mathrm{Hess} f\left(\frac{\partial}{\partial x_j},\frac{\partial}{\partial y_{a}}\right) &= \frac{1}{2} (\Lie_X g)\left(\frac{\partial}{\partial x_j},\frac{\partial}{\partial y_{a}}\right) = \frac{ \partial^2 f}{\partial x_i \partial y_a} + \frac{ \partial^2 f}{\partial y_a \partial x^i} 
    \end{align*}
    Since second partial derivatives commute, this shows that the mixed partials 
    \[ \frac{ \partial^2 f}{\partial x_i \partial y_a} = 0. \]
    This implies that $\frac{\partial f}{\partial x_i}$ is a function of $x$ only and $\frac{\partial f}{\partial y_a}$ is a function of $y$ only, so $f$ splits. 
\end{proof}

Now we consider Bach solitons on the product of a $1$-manifold and a $3$-manifold.  In a coordinate neighborhood we can then write the  metric  as $dt^2 + g_N$ where $g_N$ is a metric on a three manifold $N$ and $t \in (a,b)$.  Letting $Y,Z$ be a vectors tangent to $N$ the Bach tensor is  \cites{DK, Helliwell}
\begin{align*}
  B \left( \frac{\partial}{\partial t}, \frac{\partial}{\partial t}\right) & = -\frac{1}{12} \Delta_N S_N - \frac{1}{4} \left( |\mathrm{Ric}_N|^2 - \frac{1}{3} (S_N)^2 \right)\\ 
    B\left( \frac{\partial}{\partial t}, Y\right) &= 0 \\
    B( Y,Z) &= \frac{1}{2}(\Delta \Ric)(Y,Z) - \frac{1}{6} \mathrm{Hess} S_N(Y,Z) - 2 \Ric_N^2(Y,Z) + \frac{7}{6} S_N \Ric_N(Y,Z)\\
    & \quad \left( -\frac{1}{12} \Delta_N S_N + \frac{3}{4} |\Ric_N|^2 - \frac{5}{12} S_N^2 \right) g(Y,Z)
\end{align*}
We obtain the following structure for Bach solitons on $S^1 \times N$. 

\begin{proposition}
  If $(M,g,X)$ is a product Bach soliton  on $S^1 \times N^3$ then $\lambda \geq 0$ and (M,g) is Bach flat if $\lambda = 0$.  If $X=\nabla f$ is a gradient field then $f$ is a function on $N$ only. 
\end{proposition}
\begin{proof}
    Suppose $(M,g,X)$ is such a Bach soliton.  
    \[ \frac{1}{2}  \Lie_X g = \frac{1}{2} \left( B + \frac{1}{12} \Delta S g\right) + \lambda g. \]
    Write 
    $X = X^1 \frac{\partial }{\partial t} + X^i \frac{\partial}{\partial x_i}$ in local coordinates as above.   From the equation for the Bach tensor we have that 
    \begin{align*} 
    \frac{1}{2} (\Lie_X g)\left( \frac{\partial}{\partial t}, \frac{\partial}{\partial t}\right) & = -\frac{1}{8} \left( |\mathrm{Ric}_N|^2 - \frac{1}{3} (S_N)^2 \right) + \lambda  \\
    (\Lie_X g)\left( \frac{\partial}{\partial t}, Y\right) &= 0.
    \end{align*}
    Let  
    \[ \alpha = -\frac{1}{8} \left( |\mathrm{Ric}_N|^2 - \frac{1}{3} (S_N)^2 \right) + \lambda  \]
    and note that $\alpha= \alpha(x)$ is a function on $N$ only. 
    Moreover, 
    \[\frac{1}{2} (\Lie_X g)\left( \frac{\partial}{\partial t}, \frac{\partial}{\partial t}\right)  = \frac{\partial }{\partial t} \left(g\left(X,\frac{\partial}{\partial t}\right)\right) = \frac{\partial}{\partial t} X^1 \]
    So we have $\frac{\partial}{\partial t} X^1 = \alpha(x)$.  Integrating this around the $S^1$ factor gives that 
    \[ 0 = \int_{S^1} \frac{\partial}{\partial t} X^1 dt = \int_{S^1} \alpha(x) dt  = \alpha(x) \cdot \mathrm{length}(S^1) \]
    So we obtain that $\alpha(x) = 0$.  This gives us that

    \[ 8\lambda =  |\mathrm{Ric}_N|^2 - \frac{1}{3} (S_N)^2.  \]
By Cauchy-Schwarz the right hand side is nonnegative and equals zero if and only of $g_N$ has constant curvature, in which case $dt^2 + g_N$ is Bach flat.  

 Now if we assume in addition that $X = \nabla f$, then the argument above gives that $\frac{\partial}{\partial t} X^1 = 0$.  So $X^1$ is a function on $N$, but we also have that  $X^1 = \frac{\partial f}{\partial t}$ so $\int_{S^1} X^1= 0$ and $X^1=0$.

 Now since
 \[ 0 = \Lie_X g\left(  \frac{\partial}{\partial t} , \frac{\partial}{\partial x^i}\right) = \frac{\partial X^i}{\partial t} + \frac{\partial X^1}{\partial x_i } = \frac{\partial X^i}{\partial t} \]

 We obtain that $X^i$ is a function on $N$ only and therefore $f$ is a function on $N$ only. 

 \end{proof}

 \begin{remark}
In the gradient case, when $f$ is a function on $N$ only, we we can see from above that the the Bach soliton equation becomes
\begin{align*}
    8\lambda &=  |\mathrm{Ric}_N|^2 - \frac{1}{3} (S_N)^2 \\
    \mathrm{Hess} f &= \frac{1}{4} \Delta \Ric - \frac{1}{12} \Hess S - \Ric^2 + \frac{7}{12} S \Ric + \left( \frac{3}{8} |\Ric|^2 - \frac{5}{24} S^2 + \lambda\right) g.
\end{align*}
We do not know if there are complete, non-compact, non-constant curvature $3$-manifolds that satisfy these equations.  We do, however, observe that there are no locally homogeneous solutions.  To see this, note if there was a solution, we'd have a homogeneous solution on the universal cover.  Then by \cite{PW} $N$ would have to split as $\mathbb{R} \times L$ where $f$ is a function on the $\mathbb{R}$ factor.  But then $L$ would be a constant curvature surface.   \cite[Proposition 4.17, 4.18] {Griffin} then shows that this case can not occur. (This can also be seen directly from the equation above). 
 \end{remark}
    
Now we consider the case $\mathbb{R} \times N$ where $N$ is compact. 

\begin{proposition}  \label{Prop:4.4}
Suppose $(M,g, X)$ is a Bach soliton product on  $\mathbb{R} \times N$ where $N$ is compact.  Then $\lambda \leq 0$ with $\lambda = 0$ if and only if $g$ is Bach flat.  Moreover, if $X = \nabla f$ then $f$ is constant on $N$,  $f= 2 \lambda t^2 +  a t + b$, and $N$ has constant scalar curvature and constant norm of the Ricci tensor. 
\end{proposition}

\begin{proof}
  Tracing the Bach soliton equation gives that 
  \begin{align} \label{eqtr}  \mathrm{div} X = \frac{1}{6} \Delta S + 4 \lambda \end{align}
For each fiber $\{t\} \times N $ we can write 
  \[ \mathrm{div} X = \frac{1}{2} (\Lie_X g)\left(\frac{\partial}{\partial t},\frac{\partial}{\partial t}\right) + \mathrm{div}_N X \]
  Where by $\mathrm{div}_N X$ we mean the divergence of the restriction of $X$ to $\{t\} \times N $.  Note that $\mathrm{div}_N X$ may depend on $t$, but that 
  \[ \int_{ \{t\} \times N } \left(\mathrm{div}_N X \right)dvol_N = 0 \] for all $t$, by compactness of $N$.  Also note that $S= S_N$ is a function on $N$ so $\int_{ \{t\} \times N } (\Delta S)  dvol_N= 0$.  Plugging this into \eqref{eqtr} and integrating gives us that 
  \[ \int_{ \{t\} \times N } \frac{1}{2}(\Lie_X g)\left(\frac{\partial}{\partial t},\frac{\partial}{\partial t}\right) dvol_N = 4\lambda \mathrm{Vol}(g_N). \]
    So the average value of $\frac12(\Lie_X g)\left(\frac{\partial}{\partial t},\frac{\partial}{\partial t}\right)$ on the fiber $\{t\} \times N$ is $4\lambda$.  Then there is at least one point $(t,x)$ such that  \[\frac12(\Lie_X g)\left(\frac{\partial}{\partial t},\frac{\partial}{\partial t}\right) = 4 \lambda. \]
 Plugging this into the soliton equation gives us at this point that 
  \begin{align*}
    4 \lambda = \frac12(\Lie_X g)\left(\frac{\partial}{\partial t},\frac{\partial}{\partial t}\right) = -\frac{1}{8} \left( |\mathrm{Ric}_N|^2 - \frac{1}{3} (S_N)^2 \right) + \lambda  
  \end{align*} So that 
  \begin{align}
       \lambda = - \frac{1}{24} \left( |\Ric|^2 - \frac{S^2}{3} \right).  \label{lambdaProp5.3}
     \end{align}
  This implies that $\lambda \leq 0$.  
  
  If $\lambda = 0$ then we have 
  \[ \frac12(\Lie_X g)\left(\frac{\partial}{\partial t},\frac{\partial}{\partial t}\right)  = -\frac{1}{8} \left( |\mathrm{Ric}_N|^2 - \frac{1}{3} (S_N)^2 \right) \leq 0\] and that \[ \int_{ \{t\} \times N } \frac{1}{2}(\Lie_X g)\left(\frac{\partial}{\partial t},\frac{\partial}{\partial t}\right) dvol_N =0.  \] Which implies that $\frac{1}{2}(\Lie_X g)\left(\frac{\partial}{\partial t},\frac{\partial}{\partial t}\right) = 0$ so that $|\mathrm{Ric}_N|^2 - \frac{1}{3} (S_N)^2  = 0 $ and $g_N$ is constant curvature and $dt^2 + g_N$ is Bach flat. 

  Now in the gradient case, by Lemma \ref{Lem:Splitting}, $f = f_1(t) + f_2(x)$ and  
  \[\frac12(\Lie_X g)\left(\frac{\partial}{\partial t},\frac{\partial}{\partial t}\right) = \mathrm{Hess} f\left(\frac{\partial}{\partial t},\frac{\partial}{\partial t}\right) = \frac{d^2f_1}{dt^2}.   \]
  Since this is a function of $t$ only, the average value over the fiber $\{t\} \times N$ is $\frac{d^2f_1}{dt^2}$, so the argument above gives that $\frac{d^2f_1}{dt^2} = 4\lambda$ so \[ f_1 = 2\lambda t^2 + bt + a.  \]
  Plugging back into \eqref{eqtr} then gives
  \[\Delta( f_2 - \frac{1}{6} S_N) = 0 . \] As $N$ is compact this implies that $f_2 = \frac{1}{6} S_N + c $ for some constant $c$.  Then the Bach soliton equation on $N$  becomes
  \begin{align*}
   \frac{1}{6} \mathrm{Hess} S &= \mathrm{Hess} f = \frac{1}{4} \Delta \mathrm{Ric} - \frac{1}{12} \mathrm{Hess} S - \Ric^2 + \frac{7}{12} S \Ric \\ & + \quad \left( \lambda -\frac{1}{12} \Delta S + \frac{3}{8} |\Ric|^2 - \frac{5}{24} S^2 \right) g
  \end{align*}
  Tracing this equation gives us that 
  \begin{align*}
      \frac{1}{8} \Delta S =  \frac{1}{8} |\Ric|^2 - \frac{1}{24} S^2 + 3 \lambda
  \end{align*}
  Plugging in \eqref{lambdaProp5.3} then gives that $\Delta S = 0$.  As $N$ is compact, this implies that $S$ is constant.  Then $|\Ric|^2$ is also constant by (\ref{lambdaProp5.3}).

\end{proof}

\begin{remark} From the proof of  Proposition \ref{Prop:4.4} , we can see that if a compact $3$-manifold has constant scalar curvature and constant $|\Ric|^2$, then there is a Bach soliton on $\mathbb{R} \times N$ if the Ricci tensor satisfies the equation
\begin{align*}
    \frac{1}{4} \Delta \Ric - \Ric^2 + \frac{7}{12} S \Ric = - \frac{1}{3} \left(  |\Ric|^2 - \frac{7}{12} S^2 \right)g.
\end{align*}
Note that constant curvature metrics will always satisfy this equation.  The first author's soliton on $\mathbb{R}\times SU(2)$ is also a solution that does not have constant curvature.  She also shows that these are the only solutions to this equation on a locally homogeneous compact three manifold. 

\end{remark}

We next consider the product of two surfaces, 
$(K \times L, g_K + g_L)$ where $K$ and $L$ are each surfaces.   The formula for the Bach tensor is  \cites{DK, Helliwell}
\begin{align*}
    B (Z,W) &= \frac{1}{3} \mathrm{Hess} S_K(X,Y) - \frac{1}{3}\left( \Delta_K S_K - \frac{1}{2} \Delta_L S_L + \frac{1}{4}\left( S_K^2 - S_L^2\right) \right)g_K(Z,W) \\ 
    B(Z,U) &= 0 \\
     B (U,V) &= \frac{1}{3} \mathrm{Hess} S_L(U,V) - \frac{1}{3}\left( \Delta_L S_L - \frac{1}{2} \Delta_K S_K + \frac{1}{4}\left( S_L^2 - S_K^2\right) \right)g_L(X,Y)
\end{align*}
where $Z,W$ are tangent vectors to $K$ and $U,V$ are tangent vectors to $L$. This gives us the following system of equations for Bach solitons. 

\begin{align*}
   \frac{1}{2} \Lie_X g (Z,W) &= \frac{1}{6} \mathrm{Hess} S_K(X,Y) + \left(\lambda  - \frac{1}{8}\left( \Delta_K S_K + \frac{1}{3} S_K^2  - \left(\Delta_L S_L + \frac{1}{3} S_L^2 \right)\right)\right) g_K(Z,W) \\ 
    \frac{1}{2} \Lie_X g (Z,U) &= 0 \\
    \frac{1}{2} \Lie_X g(U,V) &= \frac{1}{6} \mathrm{Hess} S_L(U,V) + \left(\lambda  - \frac{1}{8}\left( \Delta_L S_L + \frac{1}{3} S_L^2  - \left(\Delta_K S_K + \frac{1}{3} S_K^2 \right)\right)\right) g_K(U,V) \\
\end{align*}

From these equations we can obtain the following. 

\begin{proposition} \label{Prop:ProductSurfaces} If $(M,g,X)$ is a split Bach soliton on $K \times L$, a product of surfaces,  where $K$ is compact, then  $(K,g_K)$ is a space of constant curvature.  If $X = \nabla f$ then $f$ is a function on $L$ only and either 
\begin{enumerate}
    \item  $\lambda>  0$ and $L$ is flat, or
    \item  $\lambda = 0$,  $f = \frac{1}{6}S_L + c $,  and $\Delta_L S_L + \frac{1}{3} S_L^2 = c_L$ for some $c_L>0$.
\end{enumerate}   
\end{proposition}

\begin{proof}
As $X$ splits, we know that the terms in the first equation above that depends on $L$ but not $K$ must be constant and similarly for the third equation with $L$ and $K$ switched.  This gives us that $\Delta_L S_L + \frac{1}{3} S_L^2 = c_L$ and $\Delta_K S_K + \frac{1}{3} S_K^2 = c_K$ for some constants $c_K, c_L$.  If $K$ is compact this implies that $(K,g_K)$ is constant curvature by  Lemma \ref{Lem:constantcurvature} below. 

Now if $X = \nabla f$, by Lemma \ref{Lem:Splitting} write $f= f_K + f_L$.  Since $K$ is constant curvature the equations become
\begin{align*}
\mathrm{Hess} f_K  &= \left(\lambda  - \frac{1}{8}\left( c_K - c_L \right)\right) g_K\\ 
\mathrm{Hess} f_L  &= \frac{1}{6} \mathrm{Hess} S_L(U,V) + \left(\lambda  - \frac{1}{8}\left( c_L  - c_K\right)\right) g_L \\
\end{align*}
The first equation tells us that $f_K$ is constant and that $\lambda = \frac{1}{8} \left( c_K - c_L \right)$. 
From the second equation then we obtain that 

\[ \mathrm{Hess}\left( f - \frac{1}{6} S_L \right) = 2 \lambda g_L.  \]
If $\lambda \neq 0$ we have a function $h$ such that $\mathrm{Hess} h = cg$ for a constant $c \neq 0$. This can happen on a complete surface only if $L$ is a flat metric on $\mathbb{R}^2$.  Then $c_L= 0$ and $c_K = \frac{1}{3} S_K^2$, so $\lambda >0$. 

If $\lambda = 0$ then we have $\mathrm{Hess} h = 0$, so the manifold either splits as a cylinder or $h$ is constant.  If the metric is a cylinder, $S_L=0$,and $\lambda = 0$ then $S_k= 0$ as well and the metric is flat.  Otherwise we have $f_L = \frac{1}{6} S_L + c$ and the soliton equation will be satisfied as long as $c_L = c_K$. 
\end{proof}

\begin{lemma} \label{Lem:constantcurvature}
    If $(K,g_K) $ is a compact surface and there is a constant $c$ such that  \[ c = \Delta S + \frac{1}{3} S^2\] then $S$ must be constant. 
\end{lemma}

\begin{proof}
    On a surface $\Ric = \frac{1}{2} S g$, so the Bochner identity, $\mathrm{div}(\nabla \nabla h) = \Ric(\nabla h) + \nabla \Delta h$, when applied to $h = S$, becomes 
    \begin{align*}
       \mathrm{div} (\nabla \nabla S) & = \mathrm{Ric}(\nabla S) + \nabla \Delta S  \\
       &= \frac{1}{2} S \nabla S + \nabla \Delta S \\
       &= \frac{1}{4} \nabla (S^2) + \nabla \Delta S.
    \end{align*} 
    On the other hand, \[ c = \Delta S + \frac{1}{3} S^2\] implies that $ \nabla (S^2) = -3 \nabla \Delta S$.  Combining these two equations  gives us that 
    \[ \mathrm{div} (\nabla \nabla S) = \frac{1}{4} \nabla \Delta S.\]
    On the other hand, the identity $\int_M ||\Lie_X g||^2 dvol_g = -2 \int_M \mathrm{div}(\Lie_X g)(X) dvol_g$ applied to $X = \frac{1}{2} \nabla S$ gives us that 
    \begin{align*} 
    \int_M || \mathrm{Hess} S||^2 dvol_g &= -\int_M \mathrm{div}(\mathrm{Hess} S)(\nabla S) dvol_g \\
    &= -\frac{1}{4} \int_M d(\Delta S)(\nabla S) dvol_g\\
    &= + \frac{1}{4} \int_M (\Delta S)^2 dvol_g\\
    \end{align*}
    But  by Cauchy-Schwarz, $|| \mathrm{Hess} S||^2 \geq \frac{(\Delta S)^2}{2}$.  So this identiy is only possible when $\Delta S = 0$.  Since the manifold is compact, this implies that $S$ is constant. 
\end{proof}

\begin{remark}
    From the proof, we can see that there is a solution in case (2) of Proposition \ref{Prop:ProductSurfaces} if and only if $(L,g_L)$ is a surface with $\Delta_L S_L + \frac{1}{3} S_L^2 = c_L$ for some $c_L>0$.  Clearly, there are local incomplete solutions to this equation.  For example, for a rotationally symmetric metric $dt^2 + \rho^2(t) d\theta^2$, the condition will be a fourth-order ODE in $\rho$.  On the other hand, the argument above shows there are no compact solutions. It would be interesting to know if there are complete surfaces that satisfy this equation. 
\end{remark}

Finally we can also consider extended Bach soliton of the form  $(K \times L, g_K + g_L)$. 

\[ \frac{1}{2} \Lie_X g = \frac{1}{2} B + \varphi g \]

 We have the following system of equations for extended Bach solitons. 

\begin{align*}
   \frac{1}{2} \Lie_X g (Z,W) &= \frac{1}{6} \mathrm{Hess} S_K(X,Y) + \left(\varphi  - \frac{1}{6}\left( \Delta_K S_K - \frac{1}{2} \Delta_L S_L + \frac{1}{4}\left( S_K^2 - S_L^2\right) \right)\right) g_K(Z,W) \\ 
    \frac{1}{2} \Lie_X g (Z,U) &= 0 \\
    \frac{1}{2}\Lie_X g(U,V) &= \frac{1}{6} \mathrm{Hess} S_L(U,V)+ \left( \varphi  - \frac{1}{6}\left( \Delta_L S_L - \frac{1}{2} \Delta_K S_K + \frac{1}{4}\left( S_L^2 - S_K^2\right) \right)\right)g_L(X,Y)
\end{align*}

Let $C = X - \frac{1}{3} \nabla^M S_M -\frac{1}{3} \nabla^N S_N $ then we have that 
\[ \frac{1}{2} \Lie_C g = \rho_1 g_M + \rho_2 g_N\]
where $\rho_i: K \times L \to \mathbb{R}$.  In particular,  if we pick a point $(x,y)$ and let  $\pi_K$ and $\pi_L$ be the projections from $K \times L$ onto the fibers $K \times \{y\}$ and $\{x\} \times L$ respectively, we obtain that $d\pi_K(C)$ is a conformal vector field on $K$ and $d\pi_L(C)$ is a conformal vector field on $L$.  Note, if the vector field $X$ is not split then these projections will depend on the choice of $(x,y)$, however if we know that one of these projections onto a compact factor is zero, we get the following rigidity. 

\begin{lemma}  \label{Lem:CompactSurfaceFactor} Let $(K\times L, g_K + g_L)$ be a extended Bach soliton with $K$ compact. If $d \pi_K(C) = 0$ for all $y \in L$, then $g_K$ is constant curvature, $X= X_L$ is a vector field on $L$, 
\[ \varphi = - \frac{1}{6}\left( \Delta_L S_L + \frac{1}{2} S^2_L\right) + \frac{1}{12} S_K^2 \]
a function on $L$ only, and
\[ \frac{1}{2} \Lie_X g_L = \frac{1}{6} \mathrm{Hess} S_L + \left( \frac{1}{3} S_K^2 - \Delta_L S_L - \frac{1}{3} S_L^2 \right) g \]

\end{lemma}

\begin{proof}
    Since $d \pi_K(C) = 0$, this implies that $X$ restricted to K is $\frac{1}{3} \nabla^K S_K$. In particular, $X$ must split, $X= X_K + X_L$ by Lemma \ref{Lem:Splitting}. Consequently, 
    \[ \varphi = \frac{1}{3}\left( \Delta_K S_K - \frac{1}{2} \Delta_L S_L + \frac{1}{4}\left( S_K^2 - S_L^2\right)\right).  \]
    Plugging this into the expression on $L$ gives us that 
    \[\Lie_{C} g_L = \frac{1}{2} \left( -\Delta_L S_L - \frac{1}{3} S_L^2 + \Delta_K S_K + \frac{1}{3} S_K^2\right) g_L . \]
    This implies that $\Delta_K S_K + \frac{1}{3} S_L^2$ is constant on $K$ since the other terms in the expression depend only on $L$.  As above, by Lemma \ref{Lem:constantcurvature}, this gives us that $g_K$ is constant curvature. 
\end{proof}

Whether or not a compact surface admits a conformal field is determined by its topology.  Namely, a compact surface with $\chi(K) <0$ does not admit any conformal fields so we have the following corollary. 

\begin{corollary}
    Let $(K\times L, g_K + g_L)$ be a extended Bach soliton with $K$ compact and $\chi(K) < 0$ then $g_K$ is a space of constant curvature. 
\end{corollary}

For the other compact surfaces, examples of conformal fields are provided by the Killing fields of the constant curvature metrics. Note that being a conformal vector field is invariant under conformal change of metric, so the Killing fields from a space of constant curvature will be conformal fields for all of the metrics in the conformal class. 

Conformal fields that are Killing fields in a conformal metric are called inessential.  The conformal fields on the torus, Klein bottle, and projective plane are all inessential.  Thus,  $S^2$ is the only compact surface with essential conformal fields.  

In fact, we have a gradient conformal field on compact surface if and only if it is a rotationally symmetric sphere. Moreover, in the complete case the only other examples are rotationally symmetric metrics on $\mathbb{R}^2$ or on a  cylinder, see \cite{Tashiro} or \cite[Theorem 5.7.4]{Petersen}.   This gives us the following result in the gradient case. 

\begin{proposition}
    Let $(K\times L, g_K + g_L)$ be a complete, gradient extended Bach soliton on a  product of two surfaces, then either 
    \begin{enumerate}
        \item one of $g_K$ or $g_L$ is a rotationally symmetric metric (on $S^2$,  $\mathbb{R}^2$, or $S^1 \times \mathbb{R}^1$), or 
        \item  $f = \frac{1}{3} S + c$ and $\Delta_K S_K + \frac{1}{3} S_K^2 = \Delta_L S_L + \frac{1}{3} S_L^2$ is a constant.  
    \end{enumerate}
    Moreover, if $K$ is compact then $g_K$ is either a metric of constant curvature or a rotationally symmetric metric on the sphere. 
\end{proposition}

\begin{proof}
From discussion above, we can see that, for a gradient extended Bach soliton, each factor always supports a gradient conformal field and thus either one of the factors is rotationally symmetric and we are in case (1),  or  $C = 0$, which gives that $\nabla f = \frac{1}{3} \nabla S$.  Moreover, setting $\Lie_C g = 0$ in the equation on $K$ and $L$ gives 
\begin{align*}
 \frac{1}{6}\left( \Delta_K S_K - \frac{1}{2} \Delta_L S_L + \frac{1}{4}\left( S_K^2 - S_L^2\right) \right) = \varphi = \frac{1}{6}\left( \Delta_L S_L - \frac{1}{2} \Delta_K S_K + \frac{1}{4}\left( S_L^2 - S_K^2\right) \right).
\end{align*}
Which implies that 
\[ \Delta_K S_K + \frac{1}{3} S_K^2  = \Delta_L S_L + \frac{1}{3} S_L^2. \]
Since the left hand side depends on $K$ only and the right hand side depends on $L$ only, the quantity is constant and we are in case (2). 

Now assume that $K$ is compact.  If we are in case (1) then either $g_K$ is a rotationally symmetric metric on the sphere, or $d \pi_K(C) = 0$.  If $d \pi_K(C) = 0$, Lemma \ref{Lem:CompactSurfaceFactor} implies $g_K$ is constant curvature.   If we are in case (2), then $g_K$ is constant curvature by Lemma \ref{Lem:constantcurvature}.   
\end{proof}
We leave as a question whether there are complete extended Bach solitons on $K \times L$ where $K$ is a torus, Klein bottle, projective plane, or 2-sphere where the factors do not have constant curvature.  One approach to determining whether they exist would be to express the equations above in terms of the conformal factor with a constant curvature metric.  These become fourth order PDEs for the conformal factor. Clearly there are local solutions to these equations, so the question becomes whether there are solutions that satisfy the boundary conditions necessary to extend to solutions on the entire space.

\bibliographystyle{alpha}

\end{document}